\documentclass[final,3p,times]{elsarticle}

\usepackage{graphicx} % Required for inserting images
\usepackage{amsmath,amssymb,amsfonts,amsthm,mathtools,bm}
\usepackage{booktabs}
\usepackage{xcolor}
\biboptions{numbers,sort&compress}
\usepackage{siunitx}
\usepackage{hyperref}

\theoremstyle{remark}
\newtheorem{remark}{Remark}[section]

\journal{Applied Mathematics and Computation} 
\begin{document}
	\begin{frontmatter}
		
		\title{Runge--Kutta numerical methods \\for ruin probabilities in classical risk model}
		
		\author[inst1]{George Kanakoudis}
		\author[inst2]{Lazaros Kanellopoulos}
		
		\affiliation[inst1]{
			organization={Department of Statistics and Insurance Science, University of Piraeus},
			city={Piraeus},
			country={Greece}
		}
		
		\affiliation[inst2]{
			organization={Department of Statistics and Actuarial-Financial Mathematics, University of the Aegean},
			city={Samos},
			country={Greece}
		}
		
		\renewcommand{\thefootnote}{\alph{footnote}}
		\footnotetext[1]{Corresponding author e-mail address: gkanak@unipi.gr}
		\footnotetext[2]{E-mail address: lkanellopoulos@aegean.gr}
		\renewcommand{\thefootnote}{\arabic{footnote}}
		
		\begin{abstract}
			In this paper, we study Runge--Kutta methods for the computation of ruin probabilities in the classical risk model through the associated Volterra integro-differential equation. The proposed framework combines fourth-order one-step and two-step Runge--Kutta schemes with numerical quadrature formulas to approximate the convolution term. In particular, the convolution term is approximated using Newton--Cotes and Gaussian quadrature formulas, including Simpson's 1/3 rule and Pareto-adapted Gauss--Jacobi quadrature. An equivalent reformulation of the Volterra equation as a system of ordinary differential equations is also considered. Implementations for Gamma and Pareto claim-size distributions are developed. Numerical results are presented to illustrate the accuracy of the proposed methods.
		\end{abstract}
		
	\end{frontmatter}
	
	\section{Introduction}
	
	In actuarial sciences, the ultimate ruin probability is a key measure of the risk associated with the surplus process of an insurance company, and it has been extensively studied in various risk models. However, even for the classical risk model, the evaluation of the ruin probability is not an easy task (e.g. Rolski et al. \cite{Rolski1999}, Asmussen and Albrecher \cite{AssmussenAlbrecher2010}). \par
	The paper concerns the classical risk model with surplus process 
	\begin{equation}
		\label{classical-risk-model}
		U(t) = u + ct-\sum_{i=1}^{N(t)}X_i, \quad t\geq 0, 
	\end{equation}
	where $u$ is the initial surplus, $c>0$ is the rate of premium per unit time and $S(t)=\sum_{i=1}^{N(t)}X_i$ denotes the aggregate claims up to time t with $S(t)=0$ if $N(t)=0$. More explicitly, $\{X_1, X_2,...\}$ is a sequence of independent and identically distributed (i.i.d.) non negative random variables, which represents individual claim sizes, with distribution function $P(x)=1-\overline{P}(x)=Pr(X\leq x)$ and mean $\mathbb{E}[X]=\int_{0}^{\infty}zdP(z)$. We assume that the sizes of these claims arrive at an insurer according to a Poisson process $\{N(t):t\geq 0\}$ with intensity $\lambda$. We further assume that the claims are independent of the claim-arrival process. \par
	Let $T=inf\{t: U(t)<0\}$ denote the time of ruin, then the probability of ruin can be expressed as 
	\[
	\psi(u)=Pr(T<\infty|U(0)=u), \quad u\geq 0.
	\]
	Moreover, the premium loading factor is given by $c= (1+\theta)\lambda EX$, where $\theta > 0$ is the security loading, ensuring that the model \eqref{classical-risk-model} is stable. \par
	Define also 
	\begin{equation}
		\label{non-ruin probability}
		\phi(u)=1-\psi(u)
	\end{equation}
	to be the probability that ruin never occurs starting from initial surplus $u$, a probability also known as the probability of non-ruin. In general, explicit expressions for the (non-)ruin probability are known only in some cases, for instance, when claim sizes have exponential or phase-type distribution (see Asmussen and Albrecher \cite{AssmussenAlbrecher2010}, pg 14-15). \par
	In recent years, a considerable number of contributions have focused on numerical methods for the computation of the ruin probability in the classical risk model. These papers roughly follow two main directions, namely attempting to solve the relevant integral or integro-differential equation numerically, and applying Monte Carlo estimation techniques. As for the second direction, several studies have concentrated on improving the Monte Carlo simulation techniques of ruin probability computation. Mnatsakanov et al. \cite{Mnatsakanov2008} proposed an empirical-type estimator of the Laplace transform of the ruin probability and subsequently recovered by a regularized Laplace inversion method. Within the same framework, Zhang et al. \cite{Zhang2014} studied the estimation of ruin probability in the classical ruin model, where estimator was constructed via Fourier inversion. Moreover, an estimation for ruin probability by calculating convolutions in the formula with quasi-Monte Carlo method proposed in \cite{CoulibalyLefevre2008}. For the estimation of ruin probability, many parametric, semi-paremetric and non-parametric estimators have been proposed in the literature; for references, see for example \cite{Politis2003, AsmussenBinswanger1997, Masiello2014} among others. \par
	On the other hand, a variety of numerical and semi-analytical techniques have been proposed in the ruin-theory literature. Considering the ruin probability evaluation in a classical risk model, a product-integration method to compute ruin probabilities was proposed by Ramsay and Usabel \cite{Ramsay1997}, while a Laplace transform based representation involving exponential integrals for Pareto claim sizes was derived by Ramsay \cite{Ramsay2003}. Albrecher et al. \cite{Albrecher2010} developed an efficient quadrature along with a rational-approximation procedure, motivated by the Trefethen–Weideman–Schmelzer method \cite{Trefethen2006}. Moment-based Padé approximations of ultimate ruin probabilities via exponential mixtures are studied by Avram et. al. \cite{Avram2011}. In the work of Goffard et al. \cite{Goffard2016}, an orthogonal-polynomial expansion method was proposed to approximate the ultimate ruin probability in the compound Poisson ruin model. Constantinescu et al. \cite{Constantinescu2018} present three formulations for the ruin probabilities with gamma distributed claims. The results are obtained as inverse Laplace transforms of integro-differential equations with two involving Mittag–Leffler functions and the third based on moments of the claim size distribution. Fourier-cosine methods are introduced in \cite{Chau2015} to analyze the sensitivity of ruin probabilities to perturbations in the moments of the claim-size and claim-arrival distributions. Gzyl et al. \cite{Gzyl2013} propose a maximum entropy approach based on fractional exponential moments is employed to determine the probability of ruin. For further numerical approaches to ruin probability approximation in the classical risk model, see also \cite{Mnatsakanov2015, Santana2017, Constantinescu2019, Martire2022}.
	\par
	In the present work, we investigate numerical approximations for ruin probabilities through the associated Volterra integro-differential equation (VIDE). The methods developed here are closely related to the recent literature on numerical Runge--Kutta (RK) schemes for VIDEs (see \cite{wen2024twostep, Maleknejad2004}). In particular, we investigate high-order Runge--Kutta methods for the numerical computation of ruin probabilities in the classical risk model through two complementary approaches. The first combines the classical fourth-order one-step RK method with Simpson's $1/3$ quadrature rule, together with an equivalent ordinary differential equation formulation whenever possible. The second employs two-step Runge--Kutta (TSRK) methods coupled with Gaussian quadrature formulas for the integral term of the VIDE, while in the Pareto case we use a Gauss--Jacobi quadrature adapted to the associated heavy-tailed kernel structure. \par
	
	The remainder of this paper is organized as follows. In Section \ref{section-RK}, we present the framework of the one-step and two-step Runge–Kutta methods together with quadrature approximation of the convolution term. Section \ref{Implementation-Gamma} develops the RK4 and TSRK schemes for Gamma claim-size distributions. In section \ref{Implementation-Pareto}, we consider the Pareto case using Gauss--Jacobi quadrature formulas. Section \ref{section-numerical results} provides numerical results that indicate the accuracy and efficiency of the proposed methods. Finally, Section \ref{section-conclusions} contains concluding remarks.

	\section{Runge--Kutta-based Methods for the integro-differential equation for ruin probability}
	\label{section-RK}
	This section describes the implementation of Runge--Kutta numerical methods for the approximation of the ruin probability in the classical risk model. Specifically, we consider combinations of Runge--Kutta schemes and numerical quadrature rules in order to solve the governing integro-differential equation. \par
	Consider the equation satisfied by the survival probability $\phi$ (see \cite{Rolski1999}, Section 5.3):
	\begin{equation}
		\phi'(u) = \frac{\lambda}{c}\phi(u) - \frac{\lambda}{c}\int_0^u \phi(u-z)\,dP(z),
		\label{eq:survival}
	\end{equation}
	where $\phi'(u)$ denotes the derivative of $\phi$ with respect to $u$. Using the relation for $\psi$ in $\eqref{non-ruin probability}$ and substituting into \eqref{eq:survival}, it follows that
	\begin{equation}
		\psi'(u) = \frac{\lambda}{c}\psi(u) - \frac{\lambda}{c}I(u) - \frac{\lambda}{c}\overline{P}(u),
		\label{eq:ruin_compact}
	\end{equation}
	where
	\[
	I(u) = \int_0^u \psi(u-z)\,dP(z)
	= \int_0^u \psi(z)\,p(u-z)\,dz
	\]
	denotes the convolution term and $p(x)$ is the density corresponding to $P(x)$. The density formulation is used to avoid the reflected Stieltjes differential notation. \par
	Hence, with the notation introduced above, equation $\eqref{eq:ruin_compact}$ may be expressed as
	\[
	\mathsf{f}\bigl(u,\psi(u),I(u)\bigr)
	= \frac{\lambda}{c}\psi(u) - \frac{\lambda}{c}I(u) - \frac{\lambda}{c}\overline{P}(u),
	\]
	which places $\psi'(u)=\mathsf{f}(u, \psi(u), I(u))$  within the general framework of Volterra integro-differential initial value problems:
	\begin{equation}
		\left\{
		\begin{aligned}
			y'(x) &= \mathsf{f}\left(x, y(x), \int_{x_0}^{x} K(x,s,y(s))\,ds\right),
			\quad x\in [x_0, T],\\
			y(x_0) &= y_0.
		\end{aligned}
		\right.
		\label{eq:standard_vide}
	\end{equation}
	where 
	\[
	\mathsf{f}:[x_0, T] \times \mathbb{R}\times \mathbb{R}, \quad K:\{(x,s):x_0\leq s \leq x \leq T\} \times \mathbb{R} \rightarrow \mathbb{R}
	\]
	are supposed to be sufficiently smooth to guarantee the existence and uniqueness of the solution (see Gripenberg et al. \cite[Chapter 3]{
Gripenberg1990}).\\
	Since the initial condition at $u=0$ is explicitly available in the classical risk model, the associated VIDE of the form $\eqref{eq:standard_vide}$ naturally suggests the application of Runge--Kutta methods together with suitable quadrature approximations for the integral term. In particular, fourth-order one-step and two-step Runge--Kutta discretizations are combined with Simpson's $1/3$ rule and Gaussian quadrature approximations, respectively.
	
	\subsection{Fourth-order Runge--Kutta (RK4) method with Simpson's 1/3 rule}
	
	It is well known that $\phi$ is a (proper) probability distribution with a mass of size $1-\psi(0)$, where
	\begin{equation}
		\label{initial condition psi(0)}
		\psi(0)=\frac{1}{1+\theta}, \quad \theta>0.
	\end{equation}
	By equation $\eqref{initial condition psi(0)}$, we set the initial condition $\psi(0)=\psi_0$. Let
	\[
	u_n=nh,
	\qquad
	\psi_n\approx \psi(u_n).
	\]
	For a stage displacement $\delta\geq 0$, define $g(z)=\psi(z)p(u_n+\delta-z)$. Thus, the history part of the convolution at the stage point $u_n+\delta$ is
	\begin{equation}
		\label{Hd}
		H_n(\delta)
		=
		\int_0^{u_n}
		g_{n,\delta}(z)\,dz,\qquad \delta \in \left\{0,\frac{h}{2},h\right\}.
	\end{equation}
	Now, using Simpson's $1/3$ rule, we approximate
	\[
	H_n(\delta)
	\approx
	\mathcal{S}^{(n)}_{1/3}[g_{n,\delta}].
	\]
	For $n=0$, we set $H_0(\delta)=0$. If $n$ is even, then
	\[
	\mathcal{S}^{(n)}_{1/3}[\phi]
	=
	\frac{h}{3}\left[\phi(u_0)+\phi(u_n)+4\sum_{\substack{1\leq j\leq n-1\\ j\ \mathrm{odd}}}\phi(u_j)+2\sum_{\substack{2\leq j\leq n-2\\ j\ \mathrm{even}}}
	\phi(u_j)\right].
	\]
	If $n$ is odd, Simpson's $1/3$ rule is applied on $[0,u_{n-1}]$ and the final interval is treated by the trapezoidal rule:
	\begin{equation}
		\mathcal{S}^{(n)}_{1/3}[\phi]=\frac{h}{3}
		\left[\phi(u_0)+\phi(u_{n-1})+4\sum_{\substack{1\leq j\leq n-2\\ j\ \mathrm{odd}}}\phi(u_j)+2\sum_{\substack{2\leq j\leq n-3\\ j\ \mathrm{even}}}\phi(u_j)\right]
		+\frac{h}{2} \left[\phi(u_{n-1})+\phi(u_n)\right].
		\label{eq:general-S13-odd}
	\end{equation}
	The first-stage convolution is $I_1=H_n(0)$, and hence
	\begin{equation}
		\label{k1-one-step-S1/3}
		k_1=\frac{\lambda}{c}[\psi_n-I_1-\overline{P}(u_n)].
	\end{equation}
	For the second stage, setting $Y_2=\psi_n+\frac{h}{2}k_1$, then the convolution at $u_n+h/2$ is approximated by
	\[
	I_2=H_n\left(\frac{h}{2}\right)+\frac{h}{4} \left[\psi_n p\left(\frac{h}{2}\right)+Y_2 p(0)\right],
	\]
	and
	\begin{equation}
		\label{k2-one-step-S1/3}
		k_2=\frac{\lambda}{c}\left[Y_2-I_2-\overline{P}\left(u_n+\frac{h}{2}\right)\right]. 
	\end{equation}
	For the third stage, setting $Y_3=\psi_n+\frac{h}{2}k_2$ and
	the corresponding convolution approximation is
	\[
	I_3=H_n\left(\frac{h}{2}\right)+\frac{h}{4} \left[\psi_n p\left(\frac{h}{2}\right)+Y_3 p(0) \right].
	\]
	Therefore,
	\begin{equation}
		\label{k3-one-step-S1/3}
		k_3=\frac{\lambda}{c} \left[Y_3-I_3-\overline{P}\left(u_n+\frac{h}{2}\right) \right].
	\end{equation}
	For the fourth stage, setting $Y_4=\psi_n+h k_3$, then the convolution at $u_n+h$ is approximated by
	\[
	I_4=H_n(h)+\frac{h}{6} \left[ \psi_n p(h)+4Y_3 p\left(\frac{h}{2}\right)+Y_4 p(0) \right].
	\]
	Hence,
	\begin{equation}
		\label{k4-one-step-S1/3}
		k_4=\frac{\lambda}{c} \left[Y_4-I_4-\overline{P}(u_n+h) \right].
	\end{equation}
	Finally, using the stage values defined in \eqref{k1-one-step-S1/3} -- \eqref{k4-one-step-S1/3}, the solution is advanced by
	\begin{equation}
		\label{psi_n+1}
		\psi_{n+1}=\psi_n+\frac{h}{6}[k_1+2k_2+2k_3+k_4].
	\end{equation}
	
	\begin{remark}
		An alternative implementation can be obtained by replacing Simpson’s $1/3$ quadrature rule in the history approximation with the composite Simpson $3/8$ rule. The Runge--Kutta stage equations remain unchanged, while the convolution terms are approximated accordingly. For values of $n$ not divisible by three, the remaining subintervals may be treated by a trapezoidal correction. The resulting scheme is analogous to the previous construction and is therefore omitted for brevity.
	\end{remark}
	\subsection{Two-Step Runge--Kutta (TSRK) methods with Gauss quadrature}
	A general $\mathsf{m}$-stage two-step Runge--Kutta discretization for \eqref{eq:ruin_compact} may be written in the form
	\begin{equation}
		\psi_{n+1}
		=
		\theta_1\psi_n+\theta_2\psi_{n-1}
		+ h\sum_{i=1}^{\mathsf{m}} v_i k_i^{[n-1]}
		+ h\sum_{i=1}^{\mathsf{m}} w_i k_i^{[n]},
		\quad n\geq 1,
		\label{eq:tsrk_ruin_update}
	\end{equation}
	where
	\[
	k_i^{[n]}
	=
	\mathsf{f}\bigl(u_n+c_i h,\Psi_i^{[n]},I_i^{[n]}\bigr),
	\quad i=1,\dots,\mathsf{m},
	\]
	and
	\begin{equation}
		\Psi_i^{[n]}
		=
		\delta_{i1}\psi_n+\delta_{i2}\psi_{n-1}
		+ h\sum_{j=1}^{\mathsf{m}} a_{ij}k_j^{[n-1]}
		+ h\sum_{j=1}^{\mathsf{m}} b_{ij}k_j^{[n]},
		\quad i=1,\dots,\mathsf{m}.
		\label{eq:tsrk_ruin_internal}
	\end{equation}
	where $\Psi_i^{[n]} \approx \psi(u_n + c_i h).$
	At the stage point \(u_n+c_i h\), the convolution term is decomposed into a history contribution over \([0,u_n]\) and a local contribution over \([u_n,u_n+c_i h]\). Accordingly, we define
	\begin{align}
		I_i^{[n]}
		&=
		I_n^{\mathrm{hist}}(u_n+c_i h)
		+ h\sum_{j=1}^{\mu_0} w_{ij}\,
		p\bigl((c_i-d_{ij})h\bigr)
		\Bigl(
		\eta_{ij1}\psi_n+\eta_{ij2}\psi_{n-1}
		\notag\\
		&\hspace{6em}
		+\, h\sum_{l=1}^{\mathsf{m}}\alpha_{ijl}k_l^{[n-1]}
		+ h\sum_{l=1}^{\mathsf{m}}\beta_{ijl}k_l^{[n]}
		\Bigr),
		\quad i=1,\dots,\mathsf{m},
		\label{eq:tsrk_ruin_integral_stage}
	\end{align}
	where \(I_n^{\mathrm{hist}}(u)\) denotes a quadrature approximation to the history integral
	\begin{equation}
		I_n^{\mathrm{hist}}(u)\approx \int_0^{u_n}\psi(z)\,p(u-z)\,dz.
		\label{eq:tsrk_ruin_history}
	\end{equation}
	In direct analogy with the general TSRK construction, one may write
	\begin{align}
		I_n^{\mathrm{hist}}(u_n+c_i h)
		&=
		h\sum_{j=1}^{\mu_1}\omega_j\,
		p\bigl((n+c_i-\xi_j)h\bigr)
		\Bigl(
		\psi_0+\frac{h}{2}\sum_{l=1}^{2\mathsf{m}+1}\widetilde{\gamma}_{jl}k_l^{[0]}
		\Bigr)
		\notag\\
		&+ h\sum_{\nu=1}^{n-1}\sum_{j=1}^{\mu_1}\omega_j\,
		p\bigl((n-\nu+c_i-\xi_j)h\bigr) \notag\\
		& \times \Bigl(
		\zeta_{j1}\psi_\nu+\zeta_{j2}\psi_{\nu-1}
		+ h\sum_{l=1}^{\mathsf{m}}\rho_{jl}k_l^{[\nu-1]}
		+ h\sum_{l=1}^{\mathsf{m}}\gamma_{jl}k_l^{[\nu]}
		\Bigr).
		\label{eq:tsrk_ruin_history_formula}
	\end{align}
	Thus \(I_n^{\mathrm{hist}}(u_n+c_i h)\) approximates the accumulated contribution from the previously computed interval \([0,u_n]\), whereas the second term in \eqref{eq:tsrk_ruin_integral_stage} represents a local quadrature correction over the current interval \([u_n,u_n+c_i h]\).
	
	Hence the resulting method advances the ruin probability by combining a two-step Runge--Kutta discretization of the differential term with a quadrature approximation of the convolution integral \(I(u)\). \par
	Since the TSRK scheme is a two-step method, the values \(\psi_0\) and \(\psi_1\) must be available before the recurrence can be started; in our implementation, \(\psi_1\) is computed from the prescribed initial value \(\psi_0\) by one step of the classical fourth-order Runge--Kutta method.
	
	\begin{remark}\label{ODE imp}
		Whenever the integro-differential equation can be expressed as a system of ODEs there is no need to approximate the convolution term by a Gauss quadrature. The TSRK method can then be applied directly to the resulting finite-dimensional ODE system.
	\end{remark}
	In the rest of this subsection we implement a Fourth-order TSRK on the ruin probability equation. In particular, we apply the TSRK method with order 4 and \( \mathsf{m}=1\)
	\label{TSRK-exp}
	Take \(\mu_0=\mu_1=2\), then the one-stage two-step Runge--Kutta method may be written as
	\begin{align*}
		\psi_{n+1}
		&=
		\theta_1\psi_n+\theta_2\psi_{n-1}
		+h v_1 k_1^{[n-1]}+h w_1 k_1^{[n]},
		\quad n\ge 1,
		\\
		k_1^{[n]}
		&=
		\frac{\lambda}{c}
		\Bigl(
		\delta_{11}\psi_n+\delta_{12}\psi_{n-1}
		+h a_{11}k_1^{[n-1]}+h b_{11}k_1^{[n]}
		\Bigr)
		-\frac{\lambda}{c}\,I_1^{[n]}
		-\frac{\lambda}{c}\,\overline{P}(u_n+c_1h),
	\end{align*}
	where $I_1^{[n]}$, $I_n^{\mathrm{hist}}(u_n+c_1h)$
	are given by \eqref{eq:tsrk_ruin_integral_stage} and \eqref{eq:tsrk_ruin_history_formula} for the specific choices of $p$, $\mathsf{m}$, $\mu_0$, $\mu_1$. The local order conditions corresponding to $d_n^{(1)}=O(h^5)$ are
	\begin{align*}
		1-\theta_1-\theta_2 &= 0,\\
		\frac{1}{\tau!}
		-\theta_2\frac{(-1)^\tau}{\tau!}
		-v_1\frac{(c_1-1)^{\tau-1}}{(\tau-1)!}
		-w_1\frac{c_1^{\,\tau-1}}{(\tau-1)!}
		&=0,
		\quad \tau=1,\dots,4.
	\end{align*}
	The stage order conditions corresponding to $d_{n,1}^{(2)}=O(h^4)$ are
	\begin{align*}
		1-\delta_{11}-\delta_{12} &= 0,\\
		\frac{c_1^{\,\tau}}{\tau!}
		-\delta_{12}\frac{(-1)^\tau}{\tau!}
		-a_{11}\frac{(c_1-1)^{\tau-1}}{(\tau-1)!}
		-b_{11}\frac{c_1^{\,\tau-1}}{(\tau-1)!}
		&=0,
	\end{align*}\\
	where $\tau=1,2,3 \, \mbox{ denotes the condition index in the order and interpolation conditions below}$.\\
	The quadrature conditions are
	\begin{align*}
		\sum_{j=1}^{2}\omega_j\xi_j^k &= \frac{1}{k+1},
		\quad k=0,1,2,3,\\
		\sum_{j=1}^{2}w_{1j}d_{1j}^k &= \frac{c_1^{\,k+1}}{k+1},
		\quad k=0,1,2,3.
	\end{align*}
	The consistency conditions for the starting and history approximations are
	\begin{align*}
		\frac{\xi_j^{\,\tau}}{\tau!}
		-\sum_{l=1}^{3}\widetilde{\gamma}_{jl}\frac{\widetilde{c}_l^{\,\tau-1}}{(\tau-1)!}
		&=0,
		\quad \tau=1,2,3,
		\\
		1-\zeta_{j1}-\zeta_{j2} &= 0,
		\\
		\frac{\xi_j^{\,\tau}}{\tau!}
		-\frac{(-1)^\tau}{\tau!}\zeta_{j2}
		-\rho_{j1}\frac{(c_1-1)^{\tau-1}}{(\tau-1)!}
		-\gamma_{j1}\frac{c_1^{\,\tau-1}}{(\tau-1)!}
		&=0,
		\quad \tau=1,2,3,
		\\
		1-\eta_{1j1}-\eta_{1j2} &= 0,
		\\
		\frac{\tau_{1j}^{\,\tau}}{\tau!}
		-\frac{(-1)^\tau}{\tau!}\eta_{1j2}
		-\alpha_{1j1}\frac{(c_1-1)^{\tau-1}}{(\tau-1)!}
		-\beta_{1j1}\frac{c_1^{\,\tau-1}}{(\tau-1)!}
		&=0,
		\quad \tau=1,2,3.
	\end{align*}
	
	Choosing
	\[
	\omega_1=\omega_2=\frac12,
	\quad
	\xi_1=\frac12-\frac{1}{2\sqrt3},
	\quad
	\xi_2=\frac12+\frac{1}{2\sqrt3},
	\]
	and
	\[
	w_{11}=\frac{c_1}{2},
	\quad
	w_{12}=\frac{c_1}{2},
	\quad
	d_{11}=\frac{c_1}{2}-\frac{c_1}{2\sqrt3},
	\quad
	d_{12}=\frac{c_1}{2}+\frac{c_1}{2\sqrt3},
	\]
	we get the remaining coefficients by solving the above algebraic system.
	
	\begin{remark}
		For the sixth-order TSRK scheme, we consider the case with two internal stages (\(\mathsf{m}=2\)) and choose \(\mu_0=\mu_1=3\). In this setting, three-point Gaussian quadrature formulas are used for both the local and history contributions of the convolution term. The corresponding coefficients are obtained by imposing the sixth-order local consistency conditions together with the associated stage-order and quadrature conditions, following the same procedure as in the fourth-order case. Similar sixth-order TSRK constructions for VIDEs were also considered in Section 6.2 of \cite{wen2024twostep}.
	\end{remark}
	
	\section{Implementation of Runge--Kutta methods for Gamma claims}
	\label{Implementation-Gamma}
	Gamma claim-size distributions have been considered in several ruin theory studies. For example, Constantinescu et al. \cite{Constantinescu2018} proposed three numerical procedures based on Mittag–Leffler functions, moments of the claim-size distribution and Laplace transform techniques for solving integro-differential equations. 
	Since Gamma distributions with  integer-valued shape parameter belong to phase-type distributions, explicit formulas for ruin probabilities are available, allowing us to assess the accuracy of our approximations (see Chapter 9 of \cite{AssmussenAlbrecher2010}). In the following, the integro-differential equation of the ruin probability with Gamma claim-sizes is approximated through the RK4 and TSRK methods developed for VIDEs (see Wen et al. \cite{wen2024twostep}). 
	\subsection{Implementation of RK with Simpson's 1/3 rule}
	Consider Gamma-distributed claim sizes with shape parameter $2$ and rate parameter $\beta>0$ with right tail 
	\begin{equation}
		\label{Gamma-tail}
		\overline{P}(u)=e^{-\beta u}(1+\beta u), \quad \beta>0. 
	\end{equation}
	Since the mean claim size is $2/\beta$, the safety loading relation implies $$\frac{\lambda}{c}=\frac{\beta}{2(1+\theta)}.$$ 
	In this case, the stage-dependent convolution also requires the auxiliary quantity
	\begin{equation}
		J_n \approx \int_0^{u_n} q_n(z)\,dz,
		\label{eq:gamma-q}
	\end{equation}
	where $q_n(z)=\psi(z)\beta^2 e^{-\beta(u_n-z)}$. Since, for an intermediate Runge--Kutta stage point $u_n+\delta$, where $\delta=h/2$ or $\delta=h$, the history part of the convolution is
	\begin{equation}
		H_n(\delta)=\int_0^{u_n}
		\psi(z)\beta^2(u_n+\delta-z)
		e^{-\beta(u_n+\delta-z)}
		\,dz.
		\label{RK4-Gamma-1}
	\end{equation}
	Also, equation $\eqref{RK4-Gamma-1}$ can be rewritten as
	\begin{align}
		H_n(\delta)
		&=
		\int_0^{u_n}
		\psi(z)\beta^2\bigl((u_n-z)+\delta\bigr)
		e^{-\beta(u_n-z)}e^{-\beta\delta}
		\,dz
		\notag\\
		&=
		e^{-\beta\delta}
		\left[
		\int_0^{u_n}
		\psi(z)\beta^2(u_n-z)e^{-\beta(u_n-z)}
		\,dz
		+
		\delta
		\int_0^{u_n}
		\psi(z)\beta^2 e^{-\beta(u_n-z)}
		\,dz
		\right].
	\end{align}
	The first integral is precisely the history convolution $I_n$. The second integral motivates the definition of the auxiliary quantity $J_n$ (see Equation $\eqref{eq:gamma-q}$). \par
	Consequently, the shifted history contribution can be written as
	\[
	H_n(\delta)
	=
	e^{-\beta\delta}
	\left(
	I_n+\delta J_n
	\right).
	\]
	For $n=0$, we set $I_0=J_0=0$. For $n\geq 1$, both $I_n$ and $J_n$ are computed by Simpson's $1/3$ history quadrature. Namely,
	\[
	I_n=\mathcal{S}^{(n)}_{1/3}[g_n],
	\qquad
	J_n=\mathcal{S}^{(n)}_{1/3}[q_n],
	\]
	where, if $n$ is even,
	\[
	\mathcal{S}^{(n)}_{1/3}[\phi]=\frac{h}{3}
	\left[\phi(u_0)+\phi(u_n)+4\sum_{\substack{1\leq j\leq n-1\\ j\ \mathrm{odd}}}\phi(u_j)+      2\sum_{\substack{2\leq j\leq n-2\\ j\ \mathrm{even}}} \phi(u_j)  \right].
	\]
	If $n$ is odd, Simpson's $1/3$ rule is applied on $[0,u_{n-1}]$ and the final interval is treated by the trapezoidal rule:
	\[
	\mathcal{S}^{(n)}_{1/3}[\phi]=\frac{h}{3}
	\left[\phi(u_0)+\phi(u_{n-1})+     4\sum_{\substack{1\leq j\leq n-2\\ j\ \mathrm{odd}}}\phi(u_j)+        2\sum_{\substack{2\leq j\leq n-3\\ j\ \mathrm{even}}}
	\phi(u_j)\right]
	+\frac{h}{2}
	[\phi(u_{n-1})+\phi(u_n)].
	\]
	The Runge--Kutta stages are then evaluated as follows. First, we define
	\begin{equation}
		\label{k1-one-step-Gamma}
		k_1=\frac{\lambda}{c} \left[ \psi_n-I_n-e^{-\beta u_n}(1+\beta u_n) \right].
	\end{equation}
	Now, setting $Y_2=\psi_n+\frac{h}{2}k_1$, then the convolution at $u_n+h/2$ is approximated by
	\[
	I_2=e^{-\beta h/2}\left(I_n+\frac{h}{2}J_n \right)+\frac{\beta^2h^2}{8}\psi_n e^{-\beta h/2}.
	\]
	Hence,
	\begin{equation}
		\label{k2-one-step-Gamma}
		k_2=\frac{\lambda}{c} \left[Y_2-I_2-e^{-\beta(u_n+h/2)} \left(1+\beta(u_n+h/2)\right) \right].
	\end{equation}
	For the third stage, we set
	\[
	Y_3=\psi_n+\frac{h}{2}k_2,
	\]
	and
	\[
	I_3=e^{-\beta h/2}
	\left(I_n+\frac{h}{2}J_n \right)+
	\frac{\beta^2h^2}{8}
	\psi_n e^{-\beta h/2}.
	\]
	Thus,
	\begin{equation}
		\label{k3-one-step-Gamma}
		k_3=\frac{\lambda}{c} \left[Y_3-I_3-e^{-\beta(u_n+h/2)} \left(1+\beta(u_n+h/2)\right) \right].
	\end{equation}
	Finally, setting $Y_4=\psi_n+h k_3$, then the convolution at $u_n+h$ is approximated by
	\[
	I_4=e^{-\beta h}\left( I_n+hJ_n \right)+\frac{h}{6} \left[\psi_n\beta^2 h e^{-\beta h}+4Y_3\beta^2\frac{h}{2}e^{-\beta h/2}\right],
	\]
	and the fourth stage is
	\begin{equation}
		\label{k4-one-step-Gamma}
		k_4=\frac{\lambda}{c}[Y_4-I_4-e^{-\beta(u_n+h)} \left(1+\beta(u_n+h)\right) ].
	\end{equation}
	The numerical solution is advanced by the classical fourth-order Runge--Kutta formula \eqref{psi_n+1}, using the stage values defined in \eqref{k1-one-step-Gamma}--\eqref{k4-one-step-Gamma}.

	\subsection{Implementation of TSRK method for Gamma$(2,\beta)$ Claims}
	Now, we employ the TSRK method with order 4, $\mathsf{m}=1$ according to Remark \ref{ODE imp} for Gamma-distributed claims with right tail given by $\eqref{Gamma-tail}$. In this case, the convolution term can be represented by introducing the auxiliary variable
	\[
	K_1(u)= \int_0^u \psi(z)e^{-\beta(u-z)}\,dz.
	\]
	Consequently, the scalar integro-differential equation is transformed into the
	following finite-dimensional system of ordinary differential equations:
	\[
	\begin{cases}
		\psi'(u)
		=
		\frac{\lambda}{c}(\psi(u)-I(u)- e^{-\beta u}(1+\beta u)), \\[0.4em]
		K_1'(u)
		=
		\psi(u)-\beta K_1(u), \\[0.4em]
		I'(u)
		=
		\beta^2 K_1(u)-\beta I(u).
	\end{cases}
	\]
	Equivalently, defining
	\[
	Y(u)
	=
	\begin{pmatrix}
		\psi(u)\\
		K_1(u)\\
		I(u)
	\end{pmatrix},
	\]
	the system can be written in the linear non-homogeneous form
	\[
	Y'(u)=MY(u)+g(u),
	\]
	where
	\[
	M=
	\begin{pmatrix}
		\frac{\lambda}{c} & 0 & -\frac{\lambda}{c}\\
		1 & -\beta & 0\\
		0 & \beta^2 & -\beta
	\end{pmatrix},
	\qquad
	g(u)
	=
	\begin{pmatrix}
		-\frac{\lambda}{c} e^{-\beta u}(1+\beta u)\\
		0\\
		0
	\end{pmatrix}.
	\]
	The initial condition is
	\[
	Y(0)
	=
	\begin{pmatrix}
		\psi(0)\\
		K_1(0)\\
		I(0)
	\end{pmatrix}
	=
	\begin{pmatrix}
		\dfrac{1}{1+\theta}\\
		0\\
		0
	\end{pmatrix}.
	\]
	The numerical method used in the implementation applies the two-step Runge--Kutta method directly to this ODE system. For a uniform mesh \(u_n=nh\), let \(Y_n\approx Y(u_n)\). A general $\mathsf{m}$-stage TSRK method has the form
	\[
	Y_{n+1}
	=
	\theta_1Y_n+\theta_2Y_{n-1}
	+h\sum_{i=1}^\mathsf{m} v_i k_i^{[n-1]}
	+h\sum_{i=1}^\mathsf{m} w_i k_i^{[n]},
	\]
	where
	\[
	k_i^{[n]}
	=
	F(u_n+c_i h,Y_i^{[n]}),
	\qquad
	F(u,Y)=MY+g(u),
	\]
	and the internal stages satisfy
	\[
	Y_i^{[n]}
	=
	\delta_{i1}Y_n+\delta_{i2}Y_{n-1}
	+h\sum_{j=1}^\mathsf{m} a_{ij}k_j^{[n-1]}
	+h\sum_{j=1}^\mathsf{m} b_{ij}k_j^{[n]}.
	\]
	Since \(F(u,Y)=MY+g(u)\) is linear in \(Y\), the stage equations are linear
	algebraic systems. Therefore, at each step, the current stage values can be
	obtained by solving a finite-dimensional linear system. \par
	For the fourth-order one-stage method ($\mathsf{m}=1$), the stage equation is
	\[
	Y_1^{[n]}
	=
	\delta_{11}Y_n+\delta_{12}Y_{n-1}
	+h a_{11}k_1^{[n-1]}
	+h b_{11}k_1^{[n]},
	\]
	with
	\[
	k_1^{[n]}=M Y_1^{[n]}+g(u_n+c_1h).
	\]
	Substituting the expression for \(k_1^{[n]}\) into the stage equation yields
	\[
	\bigl(\mathbb{I}_3-hb_{11}M\bigr)Y_1^{[n]}
	=
	\delta_{11}Y_n+\delta_{12}Y_{n-1}
	+h a_{11}k_1^{[n-1]}
	+h b_{11}g(u_n+c_1h),
	\]
	where \(\mathbb{I}_3\) is the \(3\times3\) identity matrix. Once \(Y_1^{[n]}\) is found,
	we compute
	\[
	k_1^{[n]}=M Y_1^{[n]}+g(u_n+c_1h),
	\]
	and update
	\[
	Y_{n+1}
	=
	\theta_1Y_n+\theta_2Y_{n-1}
	+h v_1 k_1^{[n-1]}
	+h w_1 k_1^{[n]}.
	\]
	
	For the sixth-order two-stage method $(\mathsf{m=2})$, the two current stages satisfy
	\[
	Y_i^{[n]}
	=
	\delta_{i1}Y_n+\delta_{i2}Y_{n-1}
	+h\sum_{j=1}^2 a_{ij}k_j^{[n-1]}
	+h\sum_{j=1}^2 b_{ij}k_j^{[n]},
	\qquad i=1,2,
	\]
	with
	\begin{equation}
		\label{k_i-TSRK}
		k_i^{[n]}=M Y_i^{[n]}+g(u_n+c_i h), \quad i=1,2.
	\end{equation}
	
	This leads to a block linear system for \(Y_1^{[n]}\) and \(Y_2^{[n]}\):
	\[
	\begin{pmatrix}
		\mathbb{I}_3-hb_{11}M & -hb_{12}M\\
		-hb_{21}M & I_3-hb_{22}M
	\end{pmatrix}
	\begin{pmatrix}
		Y_1^{[n]}\\
		Y_2^{[n]}
	\end{pmatrix}
	=
	\begin{pmatrix}
		R_1^{[n]}\\
		R_2^{[n]}
	\end{pmatrix},
	\]
	where
	\[
	R_i^{[n]}
	=
	\delta_{i1}Y_n+\delta_{i2}Y_{n-1}
	+h\sum_{j=1}^2 a_{ij}k_j^{[n-1]}
	+h\sum_{j=1}^2 b_{ij}g(u_n+c_jh),
	\qquad i=1,2.
	\]
	After solving this block system, the stage derivatives are evaluated as $\eqref{k_i-TSRK}$ and the new approximation is computed from
	\[
	Y_{n+1}=\theta_1Y_n+\theta_2Y_{n-1}+h\sum_{i=1}^2 v_i k_i^{[n-1]}+h\sum_{i=1}^2 w_i k_i^{[n]}.
	\]
	In the present ODE system formulation, the special Erlang structure of the Gamma$(2,\beta)$
	density makes it possible to represent the convolution term exactly by the
	auxiliary variables \(K_1(u)\) and \(I(u)\). Thus, the quadrature-based
	convolution approximation is replaced by an equivalent finite-dimensional ODE
	system.

	\section{Implementation of Runge--Kutta methods for Pareto claims}
	\label{Implementation-Pareto}
	Among heavy-tailed claim-size distributions, the Pareto distribution is of particular interest in actuarial science and ruin theory. Practical computation algorithms have been proposed in the literature for ruin probabilities under different versions of the Pareto distribution (see \cite{Ramsay1997}, \cite{Ramsay2003}, \cite{Albrecher2009}). In this subsection, we implement two numerical approaches based on the RK4 method combined with Simpson's $1/3$ quadrature rule and a two-step Runge--Kutta (TSRK) method combining the Gauss--Jacobi quadrature for the computation of ruin probabilities with Pareto claim sizes.

	\subsection{Implementation of RK with Simpson's 1/3 rule}
	\label{section-RK4-S1/3-Pareto}
	
	A single Pareto claim-size distribution with support on the positive real axis and mean $\mathbb{E}[X]=1$ is considered. Since Pareto distributions with integer parameter are commonly used in risk theory, the right tail is given by
	\begin{equation}
		\overline P(u)=\left(\frac{m}{u+m}\right)^{m+1}, \quad u\geq 0 \quad \text{and} \quad m=1,2, \ldots
		\label{pareto-tail}
	\end{equation}
	with corresponding density
	\begin{equation}
		\label{Pareto-density}
		p(u)=\frac{m+1}{m}\left( \frac{m}{u+m}
		\right)^{m+2}, \quad u\geq 0.
	\end{equation}
	Recall \eqref{Hd}-\eqref{eq:general-S13-odd}, then for the Pareto--Lomax density the first-stage convolution is $I_1=H_n(0)$ and
	\begin{equation}
		\label{k1-Pareto}
		k_1=\frac{\lambda}{c} \left[ \psi_n-I_1-\left(\frac{m}{u_n+m}\right)^{m+1} \right].
	\end{equation}
	For the second stage, setting $Y_2=\psi_n+\frac{h}{2}k_1$, the convolution at $u_n+h/2$ is approximated by
	\[
	I_2=H_n\left(\frac{h}{2}\right)+\frac{h}{4} \left[ \psi_n p\left(\frac{h}{2}\right)+Y_2 p(0) \right],
	\]
	where $p(0)$ and $p(h/2)$ are given by $\eqref{Pareto-density}$ at $u=0$ and $u=h/2$, respectively. Hence,
	\begin{equation}
		k_2=\frac{\lambda}{c} \left[Y_2-I_2-\left( \frac{m}{u_n+h/2+m}\right)^{m+1} \right].
	\end{equation}
	For the third stage, setting $Y_3=\psi_n+\frac{h}{2}k_2$, the corresponding convolution approximation is
	\[
	I_3=H_n\left(\frac{h}{2}\right)+\frac{h}{4} \left[ \psi_n p\left(\frac{h}{2}\right)+Y_3 p(0) \right].
	\]
	Thus,
	\begin{equation}
		k_3=\frac{\lambda}{c}\left[Y_3-I_3-\left(\frac{m}{u_n+h/2+m}\right)^{m+1} \right].
	\end{equation}
	For the fourth stage, setting $Y_4=\psi_n+h k_3$, then the convolution at $u_n+h$ is approximated by
	\[
	I_4=H_n(h)+\frac{h}{6} \left[ \psi_n p(h)+4Y_3 p\left(\frac{h}{2}\right)+Y_4 p(0) \right],
	\]
	where $p(h)$ is obtained from $\eqref{Pareto-density}$ for $u=h$. Therefore,
	\begin{equation}
		\label{k4-Pareto}
		k_4=\frac{\lambda}{c} \left[ Y_4-I_4-\left( \frac{m}{u_n+h+m}\right)^{m+1} \right].
	\end{equation}
	Finally, from the stage values \eqref{k1-Pareto}--\eqref{k4-Pareto}, the solution is obtained by the classical RK4 formula \eqref{psi_n+1}.

	\subsection{Implementation of TSRK method with order 4 and one internal stage}
	\label{section-tsrk-pareto-history-local}
	
	We now specialize the TSRK discretization to Pareto claims. The TSRK update
	\eqref{eq:tsrk_ruin_update} and internal stage \eqref{eq:tsrk_ruin_internal}
	are unchanged; only the convolution term is evaluated by a Pareto-adapted
	history--local quadrature.
	
	At the stage point
	\[
	s_1^{[n]}=u_n+c_1 h,
	\]
	the stage derivative is computed as
	\begin{equation}
		k_1^{[n]}
		=
		\frac{\lambda}{c}\Psi_1^{[n]}
		-
		\frac{\lambda}{c}I_1^{[n]}
		-
		\frac{\lambda}{c}\overline P\left(s_1^{[n]}\right),
		\label{eq:pareto-tsrk-stage-history-local}
	\end{equation}
	where
	\[
	I_1^{[n]}
	=
	\int_0^{s_1^{[n]}}
	\psi(z)p(s_1^{[n]}-z)\,dz .
	\]
	With the lag variable \(x=s_1^{[n]}-z\) and
	\(\delta_1=c_1 h\), this convolution is split as
	\begin{equation}
		I_1^{[n]}
		=
		\int_{\delta_1}^{s_1^{[n]}}
		\psi\left(s_1^{[n]}-x\right)p(x)\,dx
		+
		\int_0^{\delta_1}
		\psi\left(s_1^{[n]}-x\right)p(x)\,dx .
		\label{eq:pareto-history-local-split}
	\end{equation}
	The first integral contains only previously computed solution values, whereas
	the second is local to the current TSRK step.
	
	For Pareto claims we use the change of variables $y=m/(x+m)$ and $x=m(1-y)/y$. Thus, for \(0\leq a<b\), it follows that
	\begin{equation}
		\int_a^b f(x)p(x)\,dx=(m+1)\int_{Y_b}^{Y_a}
		f\!\left(m\frac{1-y}{y}\right)y^m\,dy,
		\label{eq:pareto-truncated-weight-integral}
	\end{equation}
	where
	\[
	Y_a=\frac{m}{a+m},\quad
	Y_b=\frac{m}{b+ m}.
	\]
	We approximate the last integral by a \(q\)-point Gaussian rule for the
	truncated weight \(y^m\) on \([Y_b,Y_a]\):
	\begin{equation}
		\int_a^b f(x)p(x)\,dx
		\simeq
		(m+1)
		\sum_{r=1}^q
		\widetilde\omega_r^{[a,b]}
		f\!\left(x_r^{[a,b]}\right),
		\label{eq:pareto-truncated-gaussian-rule}
	\end{equation}
	where
	\[
	x_r^{[a,b]}=m\frac{1-y_r^{[a,b]}}{y_r^{[a,b]}} .
	\]
	
	The nodes and weights are obtained from the moments
	\begin{equation}
		M_k^{[a,b]}=\int_{Y_b}^{Y_a} y^{k+m}\,dy
		=\frac{Y_a^{m+k+1}-Y_b^{m+k+1}}{m+k+1},
		\qquad k=0,\ldots,2q-1 .
	\end{equation}
	Let us define
	\[
	H_0^{[a,b]}=
	\left(M_{j+k}^{[a,b]}\right)_{j,k=0}^{q-1},
	\qquad
	H_1^{[a,b]}=
	\left(M_{j+k+1}^{[a,b]}\right)_{j,k=0}^{q-1}.
	\]
	Then \(y_r^{[a,b]}\) are the generalized eigenvalues of
	\[
	H_1^{[a,b]}v=y\,H_0^{[a,b]}v,
	\]
	and the weights satisfy
	\[
	\sum_{r=1}^q
	\widetilde\omega_r^{[a,b]}
	\left(y_r^{[a,b]}\right)^k=M_k^{[a,b]},
	\qquad k=0,\ldots,q-1 .
	\]
	
	Applying this rule to the local part of
	\eqref{eq:pareto-history-local-split}, with \(a=0\) and
	\(b=\delta_1\), gives
	\begin{equation}
		I_{1,\mathrm{loc}}^{[n]}\simeq(m+1)
		\sum_{r=1}^q
		\widetilde\omega_r^{[0,\delta_1]}
		\psi\!\left(s_1^{[n]}-x_r^{[0,\delta_1]}\right).
		\label{eq:pareto-tsrk-local-part}
	\end{equation}
	
	For the history part, split \([0,u_n]\) into
	\([u_\ell,u_{\ell+1}]\), \(\ell=0,\ldots,n-1\). The corresponding lag interval is
	\[
	a_{\ell 1}^{[n]}=s_1^{[n]}-u_{\ell+1},
	\qquad
	b_{\ell 1}^{[n]}=s_1^{[n]}-u_\ell .
	\]
	Hence,
	\begin{equation}
		I_{1,\mathrm{hist}}^{[n]}\simeq(m+1)
		\sum_{\ell=0}^{n-1}
		\sum_{r=1}^q
		\widetilde\omega_{\ell 1 r}^{[n]}
		\psi\!\left(s_1^{[n]}-x_{\ell 1 r}^{[n]}\right),
		\label{eq:pareto-tsrk-history-part}
	\end{equation}
	where
	\[
	x_{\ell 1 r}^{[n]}=m\frac{1-y_{\ell 1 r}^{[n]}}{y_{\ell 1 r}^{[n]}}
	\]
	and \(y_{\ell 1 r}^{[n]}\), \(\widetilde\omega_{\ell 1 r}^{[n]}\) are the nodes
	and weights of the same Gaussian rule on
	\[
	\left[
	\frac{m}{b_{\ell 1}^{[n]}+m},
	\frac{m}{a_{\ell 1}^{[n]}+m}
	\right].
	\]
	Thus,
	\begin{equation}
		I_1^{[n]}
		\simeq
		I_{1,\mathrm{hist}}^{[n]}
		+
		I_{1,\mathrm{loc}}^{[n]}.
		\label{eq:pareto-tsrk-total-convolution}
	\end{equation}
	The Pareto tail entering \eqref{eq:pareto-tsrk-stage-history-local} is
	\[
	\overline P(s_1^{[n]})
	=
	\left(
	\frac{m}{s_1^{[n]}+m}
	\right)^{m+1}.
	\]
	The values of \(\psi\) at off-grid points are evaluated by the cubic interpolation used in the TSRK construction: local points by the TSRK stage interpolant and history points from the already computed numerical solution. The resulting Pareto-weighted history--local scheme is denoted by TSRK4--PHL.
	
	\begin{remark}
		An alternative is to combine the convolution and tail terms into a single
		improper integral. Define
		\[
		\psi_*(s)=
		\begin{cases}
			\psi(s), & s\geq 0,\\
			1,       & s<0.
		\end{cases}
		\]
		Then
		\[
		\int_0^u \psi(u-z)p(z)\,dz+\overline P(u)
		=
		\int_0^\infty \psi_*(u-z)p(z)\,dz .
		\]
		Using the same Pareto transformation and a \(q\)-point Gauss--Jacobi rule on
		\([-1,1]\) with weight \((1+t)^m\), the nodes and weights are
		\[
		x_r= m\frac{1-t_r}{1+t_r},
		\qquad
		\mathcal W_r=\frac{m+1} {2^{m+1}}\omega_r^{(0,m)} .
		\]
		The corresponding TSRK stage is
		\[
		k_i^{[n]}=\frac{\lambda}{c}\Psi_i^{[n]}-
		\frac{\lambda}{c}\sum_{r=1}^q\mathcal W_r
		\psi_*\!\left(u_n+c_i h-x_r\right).
		\]
		In our computations, this fully Pareto-adapted improper quadrature did not improve on the TSRK--PHL approach.
	\end{remark}
	
	\section{Numerical results}
	\label{section-numerical results}
	This section summarizes the findings of comparing different orders of Runge--Kutta (RK) methods combined with various numerical integration schemes to calculate the probabilities of ruin under both Gamma and Pareto claim size distributions.\par
	We consider Gamma-distributed claim sizes with different parameter choices, summarized in Tables \ref{Table-Gamma} and \ref{Table-Gamma-2}. In particular, Table \ref{Table-Gamma} presents numerical approximations of the survival probability when claim sizes follow Gamma$(2, 2.4)$ with safety loading $\theta=0.2$. The numerical approximations obtained using the fourth-order Runge--Kutta method combined with Simpson’s $1/3$ quadrature rule (denoted by RK4--S$1/3$), together with 
	the two-step Runge--Kutta method of order four, indicated by TSRK4--G, respectively. We compare our results with the three infinite-series formulas for ruin probabilities proposed in \cite{Constantinescu2018}. We also include the results obtained by Martire in \cite{Martire2022}, where the ruin probability is computed through a numerical solution of the associated Volterra integral equation based on a Lipschitz-continuity approach. The notations M-1, M-2, and M-3 correspond to Methods 1–3 in \cite{Constantinescu2018}, associated with different discretization steps $\Delta$, while AM denotes the results obtained using the method in \cite{Martire2022}. \par 
	Table \ref{Table-Gamma-2} presents the corresponding numerical approximations of the ruin probability for a claim size follows Gamma$(2,1)$ with $\theta=1.5$. The same numerical methods are employed on the interval $u\in [0, 14.892]$ using the same step size $h=0.0016$. In this case, the results are compared with approximated values obtained through the scaled Laplace transform inversion method of Mnatsakanov et al. in \cite{Mnatsakanov2015}, as reported in \cite{Goffard2016}. The scaled Laplace transform technique has been presented in \cite{Mnatsakanov2013} and applied to the approximation of ruin probabilities in \cite{Mnatsakanov2015}. \par
	
	\begin{table}[!htbp]
		\centering
		\scriptsize
		\caption{Survival probability for claim sizes distributed as Gamma$(r=2, \beta=2.4)$ approximated via a fourth-order Runge--Kutta method with Simpson's $1/3$ and a two-step Runge--Kutta with Gaussian quadrature for $h=0.0016$ and $\theta=0.2$.}
		
		\begin{tabular}{c c cc c c c c c  }
			\toprule
			
			& Exact & \multicolumn{2}{c}{M-1} & M-2 & M-3 & AM &  RK4--S$1/3$ & TSRK4--G \\
			
			\cmidrule(lr){3-4}
			
			$u$ 
			& & $(\Delta=0.01)$ & $(\Delta=0.001)$ 
			&  &  &  &     \\
			
			\midrule
			
			0 & 0.16667  & 0.167 & 0.167 & 0.167 & 0.167 & 0.167 &   0.167 & 0.167 \\
			1  & 0.35167  & 0.350 & 0.352 & 0.352 & 0.352 & 0.352 &  0.352 & 0.352  \\
			2  & 0.50573  & 0.503 & 0.505 & 0.506 & 0.506 & 0.506 &  0.506 & 0.506  \\
			3  & 0.62347  & 0.620 & 0.623 & 0.623 & 0.623 & 0.623 &  0.623 & 0.623  \\
			4  & 0.71317  & 0.709 & 0.713 & 0.713 & 0.713 & 0.713 &  0.713 & 0.713  \\
			5  & 0.78150  & 0.777 & 0.781 & 0.782 & 0.782 & 0.781 &  0.782 & 0.782  \\
			6  & 0.83355  & 0.830 & 0.833 & 0.834 & 0.834 & 0.833 &   0.833 & 0.833  \\
			7  & 0.87321  & 0.870 & 0.873 & 0.873 & 0.873 & 0.873 &   0.873 & 0.873  \\
			8  & 0.90341  & 0.900 & 0.903 & 0.903 & 0.903 & 0.903 &   0.903 & 0.903  \\
			9  & 0.92642  & 0.923 & 0.925 & 0.926 & 0.926 & 0.926 &   0.926 & 0.926  \\
			10 & 0.94395  & 0.939 & 0.941 & 0.944 & 0.944 & 0.944 &   0.944 & 0.944  \\
			
			\bottomrule
		\end{tabular}
		\label{Table-Gamma}
	\end{table}

	\begin{table}[!h]
		\centering
		\scriptsize
		\caption{Ruin probability for claim sizes distributed as Gamma$(r=2, \beta=1)$ approximated via a fourth-order Runge--Kutta method with Simpson's  $1/3$ with $\theta=1.5$  and $h=0.0016$.}
		
		\begin{tabular}{c c c c c  }
			\toprule
			$u$ 
			& \multicolumn{1}{c}{Exact} 
			& \multicolumn{1}{c}{Scaled Laplace} 
			& \multicolumn{1}{c}{RK4--S$1/3$} 
			& \multicolumn{1}{c}{TSRK4--G}\\
			& \multicolumn{1}{c}{} 
			& \multicolumn{1}{c}{transform inversion} 
			& \multicolumn{1}{c}{} \\
			
			\midrule
			
			0.654427 & 0.320476925 & 0.320348 & 0.32048034343 & 0.32048034341  \\
			1.37683 & 0.241870349 & 0.241780 & 0.24179489538 & 0.24179489535 \\
			2.18027 & 0.173091763 & 0.173045 & 0.17305279492 & 0.17305279489  \\
			3.08527 & 0.117262707 & 0.117250 & 0.11728660791 & 0.11728660786  \\
			4.12126 & 0.074564180 & 0.0745721 & 0.07455299424 & 0.07455299419  \\
			5.33268 & 0.043754077 & 0.0437713 & 0.04375170606 & 0.04375170600  \\
			6.79131  & 0.022988447 & 0.023007 & 0.02298140146 & 0.02298140139  \\
			8.62459  & 0.010230474 & 0.0102453 & 0.01023311287 & 0.01023311279  \\
			11.0941  & 0.003436759 & 0.0034458 & 0.00343629088 & 0.00343629080  \\
			14.892   & 0.000642014 & 0.000645514 & 0.00064177551 & 0.00064222919  \\
			
			\bottomrule
		\end{tabular}
		\label{Table-Gamma-2}
	\end{table}
	
	Table \ref{Table-Pareto} presents numerical approximations of the ruin probability for Pareto claim sizes with parameter $m=1$ and for several values of the safety loading parameter $\theta$. The approximations are obtained using the classical fourth-order Runge--Kutta method combined with Simpson's $1/3$ quadrature rule (RK--S$1/3$), as described in Subsection $\ref{section-RK4-S1/3-Pareto}$, together with the two-step Runge--Kutta method combined with Gauss--Jacobi quadrature (TSRK4-PHL) presented in Subsection $\ref{section-tsrk-pareto-history-local}$. The computations are performed with step-size $h=0.01$. The numerical results are compared with exact known ruin probabilities for the Pareto case given in \cite{Ramsay2003} (denoted by RAM), as well as with the approximations obtained by Martire in \cite{Martire2022} (denoted by AM).

	\begin{table}[!htbp]
		\centering
		\tiny
		\setlength{\tabcolsep}{8pt}
		\caption{Ruin probability for Pareto claims with parameter $m=1$ approximated by the fourth-order Runge--Kutta method with Simpson's $1/3$  and the TSRK method, for various values of $\theta$ with $h=0.01$.}
		
		\begin{tabular}{c cccc cccc cccc}
			\toprule
			
			& \multicolumn{4}{c}{$\theta=0.10$}
			& \multicolumn{4}{c}{$\theta=0.25$}
			& \multicolumn{4}{c}{$\theta=1.00$} \\
			
			\cmidrule(lr){2-5}
			\cmidrule(lr){6-9}
			\cmidrule(lr){10-13}
			
			$u$ 
			& RAM & AM & RK4--S$1/3$ & TSRK4--PHL
			& RAM & AM & RK4--S$1/3$ & TSRK4--PHL
			& RAM & AM & RK4--S$1/3$ & TSRK4--PHL \\
			
			\midrule
			
			10  & 0.627128 & 0.627135 & 0.627106 & 0.627118 & 0.372677 & 0.372684  & 0.372666 & 0.372665 & 0.102523 & 0.102524 & 0.102521 & 0.102523 \\
			20  & 0.498142 & 0.498171 & 0.498089 & 0.498122 & 0.245260 & 0.245278 & 0.245238 & 0.245241 & 0.055049 & 0.055051 & 0.055047 & 0.055049 \\
			30  & 0.411437 & 0.411493 & 0.411351 & 0.411407 & 0.178338 & 0.178365 & 0.178306 & 0.178313 & 0.036887 & 0.036889 & 0.036884 & 0.036887 \\
			40  & 0.347893 & 0.347981 & 0.347775 & 0.347856 & 0.137559 & 0.137593 & 0.137520 & 0.137531 & 0.027509 & 0.027511 & 0.027506 & 0.027509 \\
			50  & 0.299155 & 0.299275 & 0.299005 & 0.299111 & 0.110519 & 0.110557 & 0.110473 & 0.110489 & 0.021847 & 0.021849 & 0.021843 & 0.021847 \\
			60  & 0.260646 & 0.260796 & 0.260465 & 0.260595 & 0.091524 & 0.091565 & 0.091473 & 0.091492 & 0.018080 & 0.018081 & 0.018076 & 0.018080 \\
			70  & 0.229551 & 0.229731 & 0.229343 & 0.229496 & 0.077594 & 0.077637 & 0.077538 & 0.077561 & 0.015402 & 0.015403 & 0.015397 & 0.015402 \\
			80  & 0.204018 & 0.204224 & 0.203784 & 0.203958 & 0.067029 & 0.067073 & 0.066969 & 0.066995 & 0.013404 & 0.013406 & 0.013400 & 0.013404 \\
			90  & 0.182761 & 0.182991 & 0.182503 & 0.182697 & 0.058794 & 0.058838 & 0.058730 & 0.058759 & 0.011859 & 0.011861 & 0.011855 & 0.011859 \\
			100 & 0.164860 & 0.165110 & 0.164579 & 0.164792 & 0.052227 & 0.052271 & 0.052160 & 0.052191 & 0.010630 & 0.010631 & 0.010625 & 0.010630 \\
			
			\bottomrule
		\end{tabular}
		
		\label{Table-Pareto}
	\end{table}

	\section{Conclusions}
	\label{section-conclusions}
	In this paper, one-step and two-step Runge--Kutta methods have been developed for approximating ruin probabilities through the associated integro-differential Volterra equation. The proposed framework provides a natural setting for the application of Runge--Kutta discretizations combined with suitable quadrature approximations for the convolution term. \par
	For the special Gamma$(2, \beta)$ case, the convolution term admits an equivalent finite-dimensional ODE representation, allowing the direct implementation of TSRK methods on the resulting system of ordinary equations. For Pareto claims, the convolution term was treated numerically through adapted quadrature procedures based on Gauss--Jacobi formulas.\par
	Numerical results indicate that the proposed methods provide accurate approximations for both light-tailed and heavy-tailed claim size distributions, remaining in close agreement with analytical formulas and benchmark numerical results in the literature. The comparisons presented in Tables~\ref{Table-Gamma}--\ref{Table-Pareto} show that the proposed RK and TSRK frameworks perform effectively across the considered examples. Moreover, the numerical results for Gamma and Pareto claim-size distributions suggest that the TSRK methods yield consistently improved approximations over the one-step schemes based on Simpson's rules. \par
	Possible extensions of the present work include the study of multistep Runge--Kutta methods in this context and applications to the Gerber--Shiu function and Sparre Andersen-type risk models with diffusion perturbations, where the corresponding equations can be formulated as Volterra integral equations of the second kind and related integro-differential equations.

	\section*{Data availability}
	
	No data was used for the research described in the article.

	\bibliographystyle{elsarticle-num}

\begin{thebibliography}{99}
		
		\bibitem{Rolski1999}
		T.~Rolski, H.~Schmidli, V.~Schmidt, J.~Teugels,
		\textit{Stochastic Processes for Insurance and Finance},
		Wiley Series in Probability and Statistics, 1999.
		
		\bibitem{AssmussenAlbrecher2010}
		S.~Asmussen, H.~Albrecher,
		\textit{Ruin Probabilities},
		2nd ed., Advanced Series on Statistical Science and Applied Probability, Vol.~14,
		World Scientific, 2010.
		
		\bibitem{Mnatsakanov2008}
		R.~M.~Mnatsakanov, L.~L.~Ruymgaart, F.~H.~Ruymgaart,
		Nonparametric estimation of ruin probabilities given a random sample of claims,
		\textit{Math. Methods Stat.} 17 (2008) 35--43.
		
		\bibitem{Zhang2014}
		Z.~Zhang, H.~Yang, H.~Yang,
		On a parametric estimator for ruin probability in the classical risk model,
		\textit{Scand. Actuar. J.} 2014~(4) (2014) 309--338.
		
		\bibitem{CoulibalyLefevre2008}
		I.~Coulibaly, C.~Lef\`evre,
		On a simple quasi-Monte Carlo approach for classical ultimate ruin probabilities,
		\textit{Insur.: Math. Econom.} 42~(3) (2008) 935--942.
		
		\bibitem{Politis2003}
		K.~Politis,
		Semiparametric estimation for non-ruin probabilities,
		\textit{Scand. Actuar. J.} (1) (2003) 75--96.
		
		\bibitem{AsmussenBinswanger1997}
		S.~Asmussen, K.~Binswanger,
		Simulation of ruin probabilities for subexponential claims,
		\textit{ASTIN Bull.} 27~(2) (1997) 297--318.
		
		\bibitem{Masiello2014}
		E.~Masiello,
		On semiparametric estimation of ruin probabilities in the classical risk model,
		\textit{Scand. Actuar. J.} 2014~(4) (2014) 283--308.
		
		\bibitem{Ramsay1997}
		C.~M.~Ramsay, M.~A.~Usabel,
		Calculating ruin probabilities via product integration,
		\textit{ASTIN Bull.} 27~(2) (1997) 263--271.
		
		\bibitem{Ramsay2003}
		C.~M.~Ramsay,
		A solution to the ruin problem for Pareto distributions,
		\textit{Insur.: Math. Econom.} 33~(1) (2003) 109--116.
		
		\bibitem{Albrecher2010}
		H.~Albrecher, F.~Avram, D.~Kortschak,
		On the efficient evaluation of ruin probabilities for completely monotone claim distributions,
		\textit{J. Comput. Appl. Math.} 233 (2010) 2724--2736.
		
		\bibitem{Trefethen2006}
		L.~N.~Trefethen, J.~A.~C.~Weideman, T.~Schmelzer,
		Talbot quadratures and rational approximations,
		\textit{BIT Numer. Math.} 46 (2006) 653--670.
		
		\bibitem{Avram2011}
		F.~Avram, D.~F.~Chedom, A.~Horv\'ath,
		On moments based Pad\'e approximations of ruin probabilities,
		\textit{J. Comput. Appl. Math.} 235~(10) (2011) 3215--3228.
		
		\bibitem{Goffard2016}
		P.-O.~Goffard, S.~Loisel, D.~Pommeret,
		A polynomial expansion to approximate the ultimate ruin probability in the compound Poisson ruin model,
		\textit{J. Comput. Appl. Math.} 296 (2016) 499--511.
		
		\bibitem{Constantinescu2018}
		C.~D.~Constantinescu, G.~Samorodnitsky, W.~Zhu,
		Ruin probabilities in classical risk models with gamma claims,
		\textit{Scand. Actuar. J.} (2018) 555--575.
		
		\bibitem{Chau2015}
		K.~W.~Chau, S.~C.~P.~Yam, H.~Yang,
		Fourier-cosine method for ruin probabilities,
		\textit{J. Comput. Appl. Math.} 281 (2015) 94--106.
		
		\bibitem{Gzyl2013}
		H.~Gzyl, P.~L.~Novi-Inverardi, A.~Tagliani,
		Determination of the probability of ultimate ruin probability by maximum entropy applied to fractional moments,
		\textit{Insur.: Math. Econom.} 53~(2) (2013) 457--463.
		
		\bibitem{Mnatsakanov2015}
		R.~M.~Mnatsakanov, K.~Sarkisian, A.~Hakobyan,
		Approximation of the ruin probability using the scaled Laplace transform inversion,
		\textit{Appl. Math. Comput.} 268 (2015) 717--727.
		
		\bibitem{Santana2017}
		D.~Santana, J.~Gonz\'alez-Hern\'andez, L.~Rinc\'on,
		Approximation of the ultimate ruin probability in the classical risk model using Erlang mixtures,
		\textit{Methodol. Comput. Appl. Prob.} 19~(3) (2017) 775--798.
		
		\bibitem{Constantinescu2019}
		C.~D.~Constantinescu, J.~M.~Ramirez, W.~R.~Zhu,
		An application of fractional differential equations to risk theory,
		\textit{Finance Stoch.} 23 (2019) 1001--1024.
		
		\bibitem{Martire2022}
		A.~L.~Martire,
		Volterra integral equations: An approach based on Lipschitz-continuity,
		\textit{Appl. Math. Comput.} 435 (2022) 127496.
		
		\bibitem{wen2024twostep}
		J.~Wen, C.~Huang, H.~Guan,
		Two-step Runge-Kutta methods for Volterra integro-differential equations,
		\textit{International Journal of Computer Mathematics}
		101~(1) (2024) 37--55.
		
		\bibitem{Maleknejad2004}
		K.~Maleknejad, M.~Shahrezaee,
		Using Runge-Kutta method for numerical solution of the system of Volterra integral equation,
		\textit{Appl. Math. Comput.} 149~(2) (2004) 399--410.
		
		\bibitem{Gripenberg1990}
		G.~Gripenberg, S.-O.~Londen, O.~Staffans,
		\textit{Volterra Integral and Functional Equations},
		Cambridge University Press, Cambridge, 1990.
		
		\bibitem{Albrecher2009}
		H.~Albrecher, D.~Kortschak,
		On ruin probability and aggregate claim representations for Pareto claim size distributions,
		\textit{Insur.: Math. Econom.} 45~(3) (2009) 362--373.
		
		\bibitem{Mnatsakanov2013}
		R.~M.~Mnatsakanov, K.~Sarkisian,
		A note on recovering the distributions from exponential moments,
		\textit{Appl. Math. Comput.} 219 (2013) 8730--8737.
		

		
	\end{thebibliography}

\end{document}